\documentclass[12pt]{amsart}
\usepackage{amsmath,amssymb,amsfonts,amscd,graphicx,xypic}
\usepackage{relsize}
\usepackage{exscale}

\usepackage{array,multirow}
\usepackage{color}

\usepackage{xcolor}

\newtheorem{lemma}{Lemma}[section]

\newtheorem{corollary}[lemma]{Corollary}
\theoremstyle{definition}

\theoremstyle{remark}

\numberwithin{equation}{section}
\theoremstyle{plain}
\newtheorem{thm}{Theorem}

\theoremstyle{definition}

\newcommand{\co}{\colon\thinspace}

\numberwithin{figure}{section}

\begin{document}

\title[Engel groups and universal surgery models]{Engel groups and universal surgery models}
\author{Michael Freedman and Vyacheslav Krushkal$^*$}

\address{Microsoft Station Q, University of California, Santa Barbara, CA 93106-6105, and
Department of Mathematics, University of California, Santa Barbara, CA 93106
}
\email{michaelf\char 64 microsoft.com}
\address{Department of Mathematics, University of Virginia, Charlottesville, VA 22904}
\email{krushkal\char 64 virginia.edu}

\thanks{$^*$ Supported in part by NSF grant DMS-1612159}

\begin{abstract}
We introduce a collection of $1/2$-${\pi}_1$-null $4$-dimensional surgery problems. This is an intermediate notion between the classically studied universal surgery models and the ${\pi}_1$-null kernels which are known to admit a solution in the topological category.
Using geometric applications of the group-theoretic $2$-Engel relation,  we show that the $1/2$-${\pi}_1$-null surgery problems are universal, in the sense that solving them is equivalent to establishing $4$-dimensional topological surgery for all fundamental groups.
As another application of these methods, we formulate a weaker version of the ${\pi}_1$-null disk lemma and show that it is sufficient for proofs of topological surgery and s-cobordism theorems for good groups.
\end{abstract}

\maketitle

\section{Introduction}
The $4$-dimensional topological surgery conjecture is known to hold for a class of {\em good} fundamental groups, and its validity for arbitrary groups remains an open problem. This conjecture, underlying geometric classification theory of topological $4$-manifolds, states that a degree one normal map $(M, \partial M)\longrightarrow (X, \partial X)$ from a topological $4$-manifold to a Poincar\'{e} pair, which is a ${\mathbb Z}[{\pi}_1 X]$-homology equivalence over $\partial X$, is normally bordant to a homotopy equivalence if and only if the Wall obstruction vanishes.   It was proved in the simply-connected case by the first author in \cite{F0}; the class of good groups has since been extended to include elementary amenable groups \cite{F1}, and more recently the groups of subexponential growth \cite{FT, KQ}. The surgery conjecture for arbitrary groups may be reduced to a collection of {\em universal} problems \cite{CF, F1} with free fundamental groups. 

The main question is whether a given surgery kernel may be ``decoupled'' from the fundamental group of the ambient $4$-manifold $M$. If the surgery kernel is represented by a ${\pi}_1$-null collection of immersed $S^2\vee S^2$s, where {\em ${\pi}_1$-null} means that all double point loops are contractible in $M$, then the surgery problem is known to admit a solution \cite[Chapter 6]{FQ}. We introduce an intermediate notion of {\em $1/2$-${\pi}_1$-null} surgery kernels, where one side is a capped grope (as in the traditional universal models), and the other side is a ${\pi}_1$-null sphere; more detailes are given below and in section \ref{classical universal section}. Our first result is that these problems are universal, in the sense that solving them is equivalent to establishing $4$-dimensional topological surgery for all fundamental groups.

\begin{thm} \label{Main theorem} \sl
The collection of $1/2$-${\pi}_1$-null surgery models is universal.
\end{thm}

\begin{figure}[ht]
%\centering
%\vspace{.5cm} 
\includegraphics[width=14.3cm]{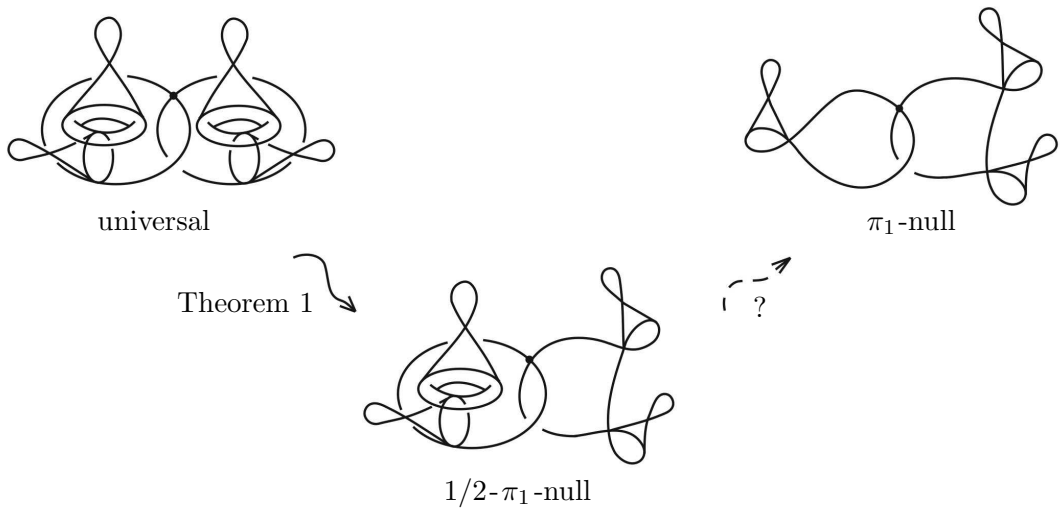}
{\small
\put(-368,94){universal}
%\put(-238,90){$1/2$-${\pi}_1$-null}
\put(-238,-8){$1/2$-${\pi}_1$-null}
\put(-78,94){${\pi}_1$-null}
\put(-338,63){Theorem \ref{Main theorem}}
\put(-120,61){?}
}
\caption{A schematic outline: Theorem \ref{Main theorem} reduces the usual universal surgery problems to 
$1/2$-${\pi}_1$-null kernels. 
In general the surfaces pictured as tori may have higher genus, and  there are higher multiplicities of double points. It is an open question whether there is a further reduction to ${\pi}_1$-null problems; a positive answer would imply the topological $4$-dimensional surgery conjecture for all groups.}
\label{Outline fig}
\end{figure}

In light of theorem \ref{Main theorem}, the current status of the surgery conjecture may be summarized as in figure \ref{Outline fig}. It is a natural question (the answer to which we do not presently know) whether the $1/2$-${\pi}_1$-null problems could be improved further, for example using the higher Engel relations discussed in \cite[Appendix]{FK}, to yield a solution to the surgery conjecture for all groups.

Next we will elaborate on the construction of the old and new universal surgery models, and the associated link-slicing problems. (See section \ref{classical universal section} for more details.)
The usual universal surgery models are given by $S^2\vee S^2$-like capped gropes \cite[Sections 2.4, 5.1]{FQ}. 
Solving these surgery problems is equivalent to the existence of free slices for links in the collection $\{$Wh(Bing(Hopf))$\}$, the links obtained from the Hopf link by (a non-trivial amount of) ramified Bing doubling and then ramified Whitehead doubling.  The free-slice condition is the requirement that the fundamental group of the complement of slices in the $4$-ball is free, generated by meridians.

The $1/2$-${\pi}_1$-null models are standard thickenings of the of $2$-complexes of the form $(G^c\vee S)\cup H^c$, where $G^c$ is a sphere-like capped grope, $S$ is a sphere with self-intersections, and $H^c$ consists of disk-like capped gropes attached to  double-point loops of $S$, figure \ref{fig:NewUniversal}. 
(Note that $S\cup H^c$ is slightly better than just a ${\pi}_1$-null sphere: contracting $H^c$ along any one of its caps gives a null-homotopy for the double point loops of $S$.)
The links corresponding to these models are of the form $\{$Wh(Bing(Wh))$\}$, ramified Whitehead double of an (iterated) ramified Bing double of the Whitehead link.

It is interesting to note that some of the representatives of $\{$Wh(Bing(Wh))$\}$ are Whitehead doubles of (homotopically trivial)$^+$ links which are known \cite{FT} to be freely slice. However, for sufficiently high genus surface stages of $G^c$ and $H^c$, and for sufficiently high multiplicity of double points, these links are of the form Wh$(L)$ where $L$ is a homotopically essential link with trivial linking numbers. In the previously considered collection $\{$Wh(Bing(Hopf))$\}$ all representatives are of this type.

The proof of theorem \ref{Main theorem} relies on a combination of geometric applications of the group theoretic $2$-Engel relation, and handle slides in a handle decomposition of surgery kernels. 
The main algebraic input in the proof is the fact that  the free group (on any number of generators) modulo the universal $2$-Engel relation $[x,[x,y]]$ is nilpotent of class $3$, see section \ref{Engel section}. This contrasts the situation with the Milnor group \cite{M}, used in the study of link homotopy: the free Milnor group on $n$ generators is nilpotent of class $n$.

The $2$-Engel relation was also used to formulate a collection of universal surgery models in \cite[Section 7]{FK} . 
Our theorem \ref{Main theorem} is a sharpening of the result of \cite{FK}, in particular the new $1/2$-${\pi}_1$-null surgery models have a very clean description in terms of their spines and the fundamental group. The proof of theorem \ref{Main theorem} is different from the argument in \cite{FK}, although both use the $2$-Engel relation.

Another application of these techniques concerns the ${\pi}_1$ null disk lemma (NDL), a technical statement underlying the proofs of the $4$D topological surgery and s-cobordism theorems for good groups.  This lemma concerns immersed disks bounded by the attaching curve of a capped grope. In section \ref{NDL section} we formulate a weaker version of NDL, where instead of the attaching curve one considers curves which are trivial in the Milnor group of the body of a capped grope.
Theorem \ref{pi1null thm} shows that this apriori weaker statement is in fact equivalent to NDL.

Section \ref{Engel section} summarizes the background material on the Milnor group and the $2$-Engel relation.  
In section \ref{classical universal section} we give a precise description of the $1/2$-${\pi}_1$-null surgery models. The proof of theorem \ref{Main theorem} is given in \ref{proof section}. Applications to the  ${\pi}_1$-null disk lemma are discussed in section \ref{NDL section}.

\section{The $2$-Engel relation and the Milnor group} \label{Engel section}

 This section gives a brief summary of the background material on link homotopy and group theory needed for the proof of theorem \ref{Main theorem}. The reader is referred to \cite{M, FK} for a detailed exposition.

Let $\pi$ be a group normally generated by elements $g_1, \ldots,  g_k$. The 
{\em Milnor group} of ${\pi}$, defined with respect to the given normal generating set $\{ g_i\}$,
is defined by
\begin{equation} \label{Milnor group definition}
M{\pi} = {\pi} \, /\, \langle\! \langle \, [g_i, g_i^y] \;\; i=1,\ldots, k,  \; \, y \in G \rangle \! \rangle.
\end{equation}

The Milnor group $ML$ of a link $L$ in $S^3$ is set to be the Milnor group $M{\pi}$ where
${\pi}={\pi}_1(S^3\smallsetminus L)$, defined with respect to meridians to the link components.

\begin{figure}[ht]
\includegraphics[height=4.5cm]{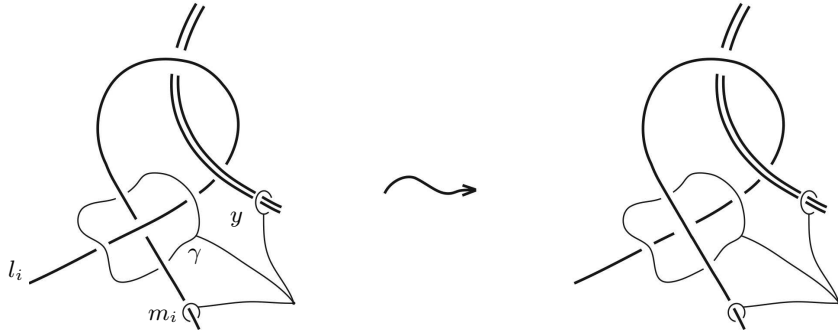}
{\scriptsize
\put(-262,6){$m_i$}
\put(-247,29){$\gamma$}
\put(-231,43){$y$}
\put(-315,22){$l_i$}
}
 \caption{A non-generic time during a link homotopy. The based curve ${\gamma}$ in the link complement corresponds to the defining relation $[m_i, m_i^y]$ of the Milnor group.}
\label{link homotopy figure}
\end{figure}

Two links are {\em link-homotopic} (figure \ref{link homotopy figure})if they
are connected by a $1$-parameter family of link maps where different components
stay disjoint for all values of the parameter.
If $L$, $L'$ are link-homotopic then 
their Milnor groups $ML$, $ML'$ are isomorphic  \cite{M}. Moreover, a $k$-component link $L$ is
link-homotopic to the $k$-component unlink  ({\em h-trivial}) if and only if $ML$ is isomorphic to the free Milnor group $MF_{m_1,\ldots, m_k}$. A $k$-component link $L$ is {\em h-trivial}$^+$  if each one of the $k$ links obtained by adding to $L$ a parallel copy of a single component is homotopically trivial.  If a link is not h-trivial, it is called {\em h-essential}.

To fix the notation, the lower central series of a group $G$ is defined inductively by $G^1=G, G^n=[G^{n-1}, G]$.
We will use will the concise notation $[g_1, g_2, \ldots, g_n]$ for the commutator $[[\ldots [g_1, g_2],\ldots, g_{n-1}], g_n]$ of elements $g_1, g_2, \ldots, g_n\in G$.

\begin{figure}[ht]
%\centering
\includegraphics[height=4.5cm]{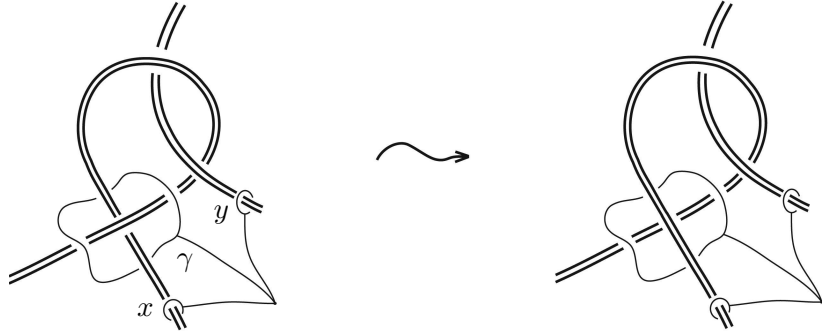}
{\small
\put(-260,7){$x$}
\put(-245,26){$\gamma$}
\put(-231,44){$y$}
}
 \caption{An elementary weak homotopy. (Compare with figure \ref{link homotopy figure}.)}
\label{weak homotopy figure}
\end{figure}

The proof of theorem \ref{Main theorem} relies on geometric applications of the $2$-Engel relation 
\begin{equation} \label{2Engel} [y,x,x]=1, {\rm or\; \; equivalently}\; \; [x,x^y]=1.
\end{equation} 
A {\em $2$-Engel group} ${\pi}$ is a group satisfying this relation for all $x,y\in {\pi}$.
Note that the defining Milnor relation (\ref{Milnor group definition}) is the same  as the $2$-Engel relation, but applied only to $x$ in a fixed set of normal generators. A geometric analogue of the $2$-Engel relation is the notion of {\em weak homotopy}, see figure \ref{weak homotopy figure} and \cite[Section 6]{FK}.

The free Milnor group on $n$ generators $MF_n$ is nilpotent of class $n$ \cite{M}. Building on earlier work of Burnside \cite{Burnside},  Hopkins \cite{Hopkins} (also see \cite{Levi}) proved that the nilpotency class of $2$-Engel groups is independent of the number of generators:

\begin{lemma}  \label{nilpotent class} \sl
Any $2$-Engel group is nilpotent of class $\leq 3$.
\end{lemma}

A proof in the setting of the Milnor group is given in \cite{FK}. The following corollary  will be used in the proof of theorem \ref{Main theorem}.

\begin{corollary} \cite[Corollary 2.3]{FK} \label{Engel corollary} \sl
Suppose $\pi$ is a group normally generated by $g_1,\ldots, g_n$. Let $g\in {\pi}^k$ be an element of the $k$-th term of the lower central series, $4\leq k\leq n$. Then $g$ may be represented in the Milnor group $M{\pi}$ as a product of (conjugates of) $k$-fold commutators $C$ of the form $[h_1,\ldots, h_k]$ where two of the elements $h_i$ are equal to each other and to a product of two generators, $h_j = h_m=g_{i_1} g_{i_2}$ for some $j\neq m$, and each other element $h_i$ is one of the generators $g_1, \ldots, g_n$.  
\end{corollary}

\begin{figure}[ht]
\centering
%\vspace{.5cm} 
\includegraphics[width=12.7cm]{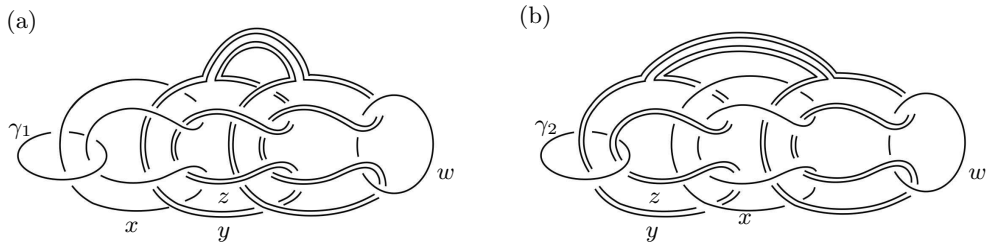}
{\scriptsize
\put(-364,33){${\gamma}_1$}
\put(-320,-4){$x$}
\put(-285,-6){$y$}
\put(-285,7){$z$}
\put(-202,16){$w$}
\put(-165,34){${\gamma}_2$}
\put(-123,-6){$y$}
\put(-1,16){$w$}
\put(-88,-1){$x$}
\put(-122,7){$z$}
\put(-364,73){(a)}
\put(-170,75){(b)}
}
\caption{(a): $\!{\gamma}_1=[x,yz,yz,w]$, $\;$ (b): $\!{\gamma}_2=[yz,x,yz,w]$.}
\label{Elementary links figure}
\end{figure}

The commutators $[h_1,\ldots, h_k]$ in the statement of the corollary have a geometric interpretation. Figure \ref{Elementary links figure} shows the {\em elementary Engel links} where the curves ${\gamma}_i$ read off these commutators in the complement of the other components, where $k=4$. These links have an important property that a handle slide involving two parallel components results in a split link consisting of the unknot and an h-trivial link. This property will be used in the proof of Theorem \ref{Main theorem}.

\section{Surgery problems} \label{classical universal section}

Our set up is a degree one normal map $(M, \partial M)\longrightarrow (X, \partial X)$ from a topological $4$-manifold to a Poincar\'{e} pair, which is a ${\mathbb Z}[{\pi}_1 X]$-homology equivalence over $\partial X$, with the vanishing Wall obstruction. A solution is a normal cobordism to a simple homotopy equivalence.

\begin{figure}
\includegraphics[height=3.7cm]{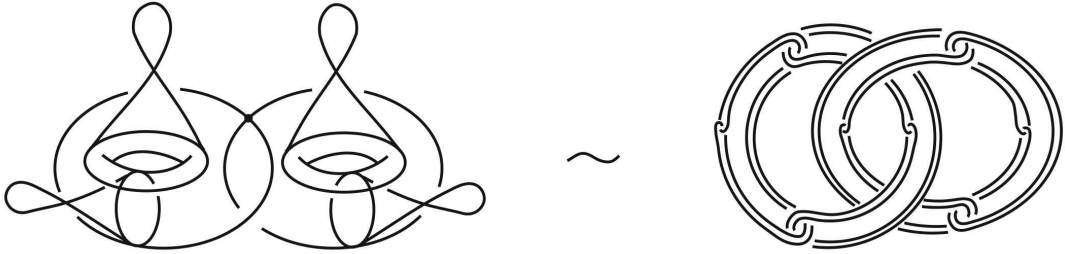} 
 \caption{The figure on the left a surgery kernel given by a height$=1$ $S^2\vee S^2$-like capped grope.  The corresponding link on the right  is obtained by Bing doubling both components of the Hopf link, and then Whitehead doubling each component.}
\label{fig:SurgeryAndGrope}
\end{figure}

In attempting to represent the surgery kernel by an embedded hyperbolic pair $S^2\vee S^2$, one finds in $M$ $2$-complexes which ``approximate'' $S^2\vee S^2$.  Specifically, one may represent \cite[Sections 2.4, 5.1]{FQ} the surgery kernel by $S^2\vee S^2$-like  capped gropes (figure \ref{fig:SurgeryAndGrope}), and the neighborhoods  of these standard examples are  themselves sources of surgery problems. This set constitutes a countable collection  of surgery problems, called {\em universal}. If  these are solvable then all unobstructed problems are solvable.

\begin{figure}[h]
\includegraphics[height=3.4cm]{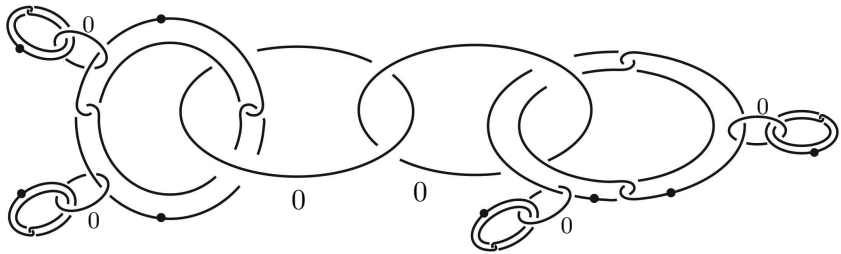}
{\small
    \put(-210,18){$0$}
    \put(-164,21){$0$}
}
{\scriptsize
\put(-289,85){$0$}
\put(-287,11){$0$}
\put(-108,9){$0$}
\put(-34,54){$0$}
}
\caption{A Kirby diagram for the surgery kernel in figure \ref{fig:SurgeryAndGrope}.}
\label{fig:wedge}
\end{figure}

There is a collection of link-slicing problems corresponding the universal surgery problems. 
Starting with $S^2\vee S^2$-like capped gropes, one obtains the collecyion $\{\!${\sl Wh(Bing(Hopf))}$\!\}$. These links are obtained from the Hopf link by a (non-zero) amount of iterated ramified Bing doubling, and then a single application of ramified $\pm$ Whitehead doubling to {\em each} component. See figure \ref{fig:SurgeryAndGrope} (right)  for an example. These links are obtained by drawing a Kirby diagram of the surgery kernel (see figure \ref{fig:wedge}), eliminating the zero-framed Hopf link in the middle (a diffeomorphism on the boundary), and canceling all $2$-handles. Solving a surgery problem produces a ``free'' slice complement, i.e. one with free ${\pi}_1$, freely generated by meridians, for {\sl Wh(Bing(Hopf))}, i.e. constructs a manifold $\simeq\vee S^1$ (generators = meridians) with boundary ${\mathcal S}_0${\sl (Wh(Bing(Hopf)))}.

\begin{figure}[ht]
%\centering
%\vspace{.5cm} 
\includegraphics[width=11.5cm]{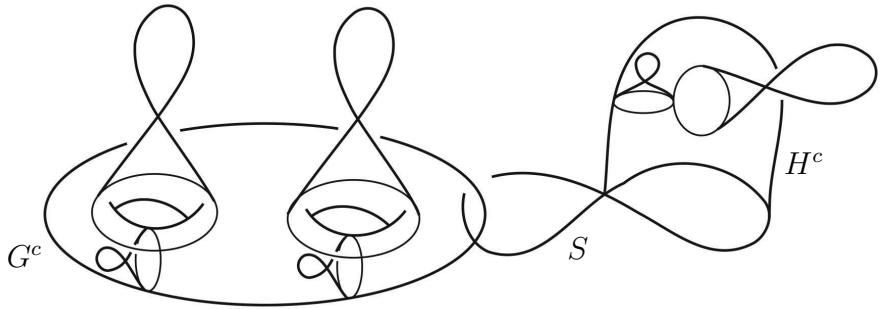}
\put(-337, 20){$G^c$}
\put(-125,23){$S$}
\put(-43,55){$H^c$}
\caption{$1/2$-${\pi}_1$-null surgery kernels of the form $(G^c\vee S)\cup H^c$. This illustration shows capped gropes $G^c, H^c$ of height $1$ (capped surfaces), and a single double point for the sphere and for the caps. In general the number of double points and the genera of surfaces (other than the sphere) are parameters in the collection of surgery models.}
\label{fig:NewUniversal}
\end{figure}

The subject of this paper, $1/2$-${\pi}_1$-null surgery kernels, are of the form $(G^c\vee S)\cup H^c$ where $G^c$ is a sphere-like capped grope, and $S$ is a sphere with self-intersections, whose double point loops bound capped gropes $H^c$, figure \ref{fig:NewUniversal}.

\begin{figure}
\includegraphics[height=3.5cm]{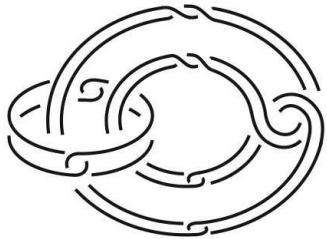} 
 \caption{The Bing double of the Whitehead link}
\label{fig:BingWh}
\end{figure}

\begin{figure}
\includegraphics[height=3cm]{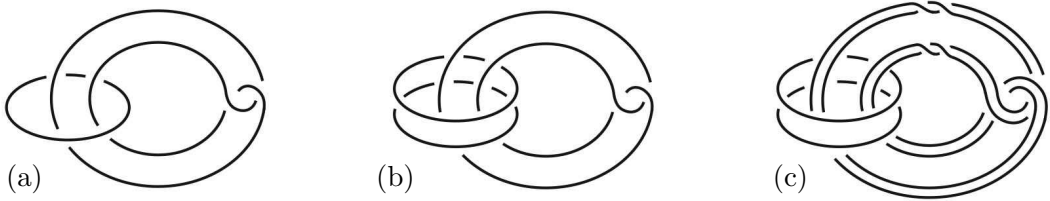} 
 \small
 \put(-400, 10){(a)}
 \put(-260,10){(b)}
 \put(-110,10){(c)}
 \caption{(a) The Whitehead link, $Wh$, is h-trivial$^+$. 
 (b) $Wh$ with a parallel copy of one of its components is h-trivial but not h-trivial$^+$. (c) $Wh$ with parallel copies of both components is h-essential.}
\label{fig:3Whiteheads}
\end{figure}

The link slicing problems, corresponding to these surgery problems, are of the form $\{$Wh(Bing(Wh))$\}$, obtained from the Whitehead link (figure \ref{fig:3Whiteheads} (a)) by first iterated ramified Bing doubling and then a single application of ramified Whitehead doubling.
These links are constructed from the Kirby diagrams of the surgery kernel analogously to figure \ref{fig:wedge}. An example of such a link (corresponding to both capped surfaces $G^c, H^c$ of genus $1$) is the Whitehead double of the link in figure \ref{fig:BingWh}. This link is freely slice \cite{FT1} since it is the Whitehead double of a homotopically trivial$^+$ link. However, links in the family $\{$Wh(Bing(Wh))$\}$ corresponding to higher genus surfaces and higher multiplicity of double points, are of the form Wh(h-essential), since ramification is used in their construction (figure \ref{fig:3Whiteheads} (b,c)).  These links are not known to be slice.

\section{Proof theorem \ref{Main theorem}} \label{proof section}

The starting point of the proof is one of the usual universal surgery problems $X$ (discussed in section \ref{classical universal section}), with a spine $A\vee B$ where $A$ and $B$ are capped gropes. The grope height raising procedure  \cite[Proposition 2.7]{FQ}  may be used to arrange $A, B$ to be capped gropes of any given height.  To be concrete, let $A$ be a $2$-stage capped grope, and $B$ a $4$-stage one. Solving this surgery problem amounts to finding a $4$-manifold $M$ which is a homotopy $1$-complex and with  with the  same boundary as $X$\footnote{A collection of ${\pi}_1$ generators of the $1$-complex is required to match up with the meridians of the link implicit in the definition of $A\vee B$ (see section \ref{classical universal section}).}.

\begin{figure}[h]
\includegraphics[height=5cm]{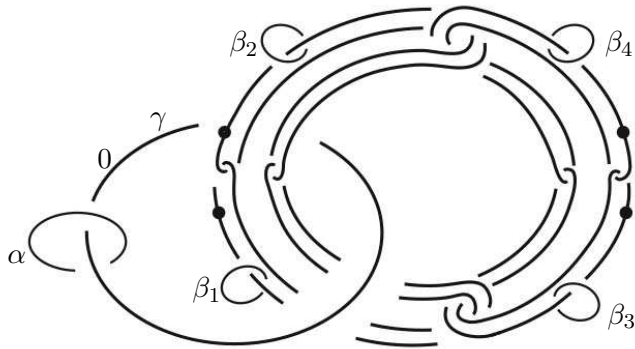}
{\small
    \put(-244,38){${\alpha}$}
\put(-174,25){${\beta}_1$}
\put(-160,118){${\beta}_2$}
\put(-17,15){${\beta}_3$}
\put(-18,118){${\beta}_4$}
\put(-210, 74){$0$}
\put(-190,90){$\gamma$}
}
\caption{The bottom two stages $H$ of the capped grope $B$.  The Kirby diagram for $H$ is given by the zero-framed curve $\gamma$ and the four dotted curves.  (The surface stages of $H$ in this example have genus $1$.) ${\beta}_i$ are the attaching curves of  $2$-stage capped gropes which form higher stages of $B$. ${\alpha}$ is the attaching curve of the $2$-stage capped grope which (capped off with the standard disk bounded by $\alpha$ in $D^4$) forms $A$.}
\label{fig:wedge1}
\end{figure}

\begin{figure}[ht]
%\centering
%\vspace{.5cm} 
\includegraphics[width=12cm]{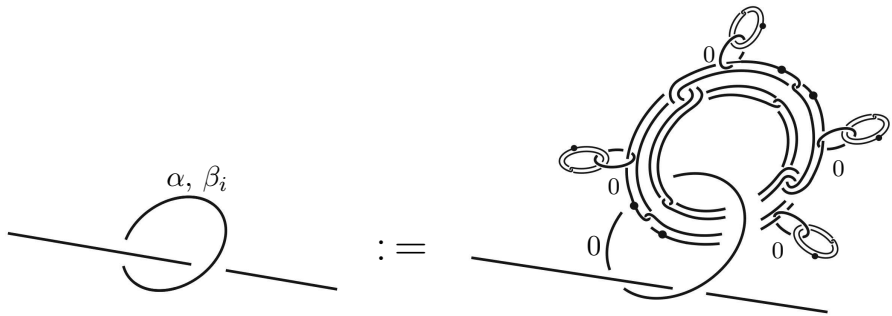}
\Large {
\put(-198, 25){$:\, =$}
}
\small{
\put(-118, 27){$0$}
\put(-277,53){${\alpha}, \, {\beta}_i$}
}
\scriptsize{
\put(-110,50){$0$}
\put(-74,100){$0$}
\put(-48,26){$0$}
\put(-25,57){$0$}
}
\caption{Explanation of the notation in figure \ref{fig:wedge1}.}
\label{fig:Notation}
\end{figure}

Consider $H$,  a disk-like symmetric grope of height $2$ consisting of the bottom two stages of $B$, figure \ref{fig:wedge1}.
To get a complete description of the surgery kernel, one replaces the curves ${\alpha}, {\beta}_i$ in figure \ref{fig:wedge1} with Kirby diagrams shown in figure \ref{fig:Notation}. These handle diagrams correspond to the case where all surface stages have genus one, and each cap has a single self-intersection point. (This case encapsulates the main features of the unresolved surgery conjecture for free groups.) In full generality, ramified Bing and Whitehead doubles are used in the construction of the dotted curves in the figures.

The attaching curve $\gamma$ of  $H$  is in the second term of the derived series ${\pi}_1(H)^{(2)}$, cf. \cite[Lemma 1.1]{FT}. 
Since the $n$-th term of the derived series of a group is contained in the $2^n$-th term of its lower central series,  ${\gamma}\in {\pi}_1 (H)^4$.

Denote by $U$ the unlink representing the $1$-handles of $H$.
In the genus $1$ case in figure \ref{fig:wedge1},  for a suitable choice of meridians $m_i$ to the dotted components $U$, the curve $\gamma$  represents the $4$-fold commutator $\big[ [m_1,m_2],[m_3,m_4] \big]$ in the  free fundamental group  ${\pi}_1(S^3\smallsetminus U)$. For brevity of notation we continue the proof in this case; it goes through for arbitrary genera of the surface stages of $H$ since the main algebraic fact, relevant in the proof, is that the attaching curve $\gamma$ is in ${\pi}_1 (H)^4$ .

\begin{figure}[ht]
%\centering
%\vspace{.5cm} 
\includegraphics[height=8.5cm]{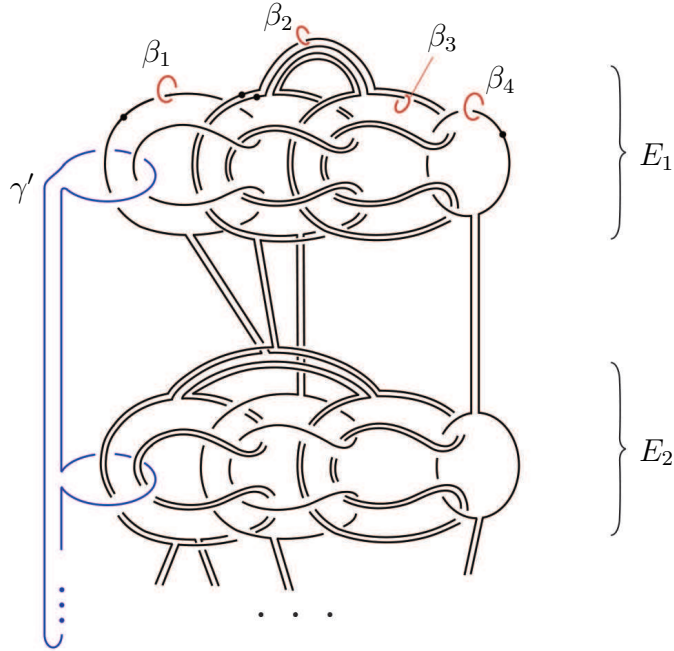}
\put(-235,175){${\gamma}'$}
\put(3,187){$E_1$}
\put(3,75){$E_{2}$}
\put(-185,226){${\beta}_{1}$}
\put(-139,240){${\beta}_{2}$}
\put(-77,232){${\beta}_{3}$}
\put(-55,216){${\beta}_{4}$}
\caption{Representation of ${\gamma}'$ as a band sum of elementary Engel links.} 
%The connecting bands are shown schematically; a precise description is obtained as the result of an isotopy taking  the link in figure \ref {fig:EngelProduct1} to the link above.}
\label{fig:EngelProduct1}
\end{figure}

By lemma \ref{Engel corollary}, $\gamma$ may be expressed in the free Milnor group $MF_{m_1,\ldots, m_4}$  as a product 
\begin{equation} \label{gamma1}
{\gamma}\, =\, \prod_k E_k^{f_k},
\end{equation}
of conjugates of commutators of the form $E_k:=[h_1,\ldots, h_4]$, where two of the elements $h_i$ are equal to each other and to a product of two generators, $h_j = h_l=m_{i_1} m_{i_2}$ for some $j\neq l$, and each other element $h_n$ is one of the generators $m_i$. Next we will construct a curve ${\gamma}'$ in the complement to the $4$-component unlink $U$, representing in $MU=MF_{m_1,\ldots, m_4}$ the same group element as $\gamma$, and whose geometry reflects the product structure in equation (\ref{gamma1}). There is a link (figure \ref{Elementary links figure}) associated to each commutator $E_k$. Consider a band sum of these elementary Engel links, one for each commutator $E_k$ in the expression (\ref{gamma1}), see figure \ref{fig:EngelProduct1}.

Note that if the curve ${\gamma}'$ is omitted in figure \ref{fig:EngelProduct1}, the remaining $4$-component link is the unlink, so putting the dots on them makes sense. Since the band-sum decomposition of $U$ will be important in the rest of the proof, we specify that the attaching curves ${\beta}_i$  of the higher stages of the capped grope B, linking the components of $U$, are positioned in the first constituent link (corresponding to $E_1$), as shown in figure \ref{fig:EngelProduct1}.

\begin{figure}[ht]
%\centering
%\vspace{.5cm} 
\includegraphics[height=6.2cm]{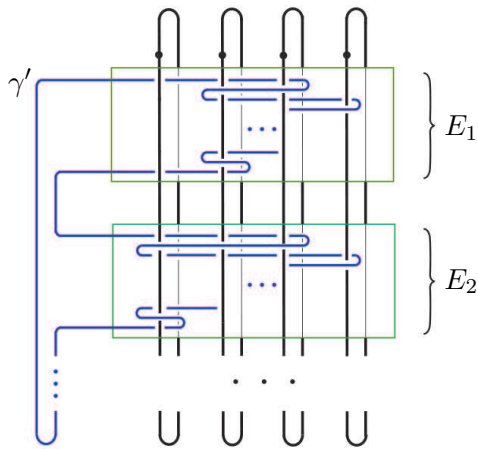}
\small
\put(-163,140){${\gamma}'$}
%\put(-100,178){${\beta}_{1}$}
%\put(-75,178){${\beta}_{2}$}
%\put(-50,178){${\beta}_{3}$}
%\put(-25,178){${\beta}_{4}$}
\put(2,123){$E_1$}
\put(2,65){$E_{2}$}
\caption{${\gamma}'$ reads off the product in equation (\ref{gamma1}) in the free Milnor group $MF_{m_1,\ldots, m_4}$, where $m_i$ are meridians to the dotted components.} %(This is a rough illustration of ${\gamma}'$, focusing on its essential features. The conjugating elements $f_i$ in the expression (\ref{gamma}) are not implemented in the figure.)}
\label{fig:EngelProduct}
\end{figure}

Straightening the dotted components by an isotopy (and letting ${\gamma}'$ evolve by an isotopy in their complement), one gets a link of the type sketched in figure \ref{fig:EngelProduct}. In fact, one could define ${\gamma}'$ as a simple closed curve reading off the product of the $E_k$ as in figure \ref{fig:EngelProduct} and carry out the rest of the proof in the group-theoretic setting, where handle slides are replaced with change of generators in the free group. We present a version using elementary Engel links, where the argument is geometrically transparent.

The commutators under consideration,  ${\gamma}=\big[ [m_1,m_2],[m_3,m_4] \big]$ and all of  the $E_i$, are of maximal length in $MF_4:=MF_{m_1,\ldots, m_4}$. Indeed, $MF_4$ is nilpotent of class $4$, and all of these elements are in the $4$-th term of the lower central series, $(MF_4)^4$. It follows that for calculations taking place in the abelian group $(MF_4)^4$ the conjugation in equation (\ref{gamma1}), as well as the choice of bands in figure \ref {fig:EngelProduct1}, are immaterial (also see \cite[Remark 4.2]{FK}).\footnote{In the general case of higher genus surface stages in $H$ the conjugation is relevant; the curve ${\gamma}'$ in figure \ref {fig:EngelProduct} and the bands in figure \ref {fig:EngelProduct1} are then drawn accordingly to represent  the expression (\ref{gamma1}).}

Since ${\gamma}, {\gamma}'$ represent the same element in $MF_4$, in the free group $F_4$  $\; [{\gamma}]\cdot [{\gamma}']^{(-1)}$ equals a product of conjugates of the Milnor relations (\ref{Milnor group definition}) $W_j$ of the form
$[m_i, m_i^y]$:
\begin{equation} \label{gamma}
[{\gamma}]\, =\, [{\gamma}'] \cdot  \prod_j W_j^{h_j}.
\end{equation}

We incorporate the new factor $W:=\prod W_j^{h_j}$ by adding one more block, reading off the word $W$ in the free group,  at the bottom to ${\gamma}'$ in figures \ref {fig:EngelProduct1}, \ref {fig:EngelProduct}. (One could also give a more precise geometric description of individual elements $W_j$ by adding specific links, similar to the elementary Engel links in figure \ref{Elementary links figure}. An important difference is that the links corresponding to $W_j$ are {\em h-trivial}, while the the elementary Engel links are h-essential. This explains why the additional factor $W$ will play a very limited role in the proof below.) Abusing the notation, we denote the links corresponding to $E_k$ and $W$ by the same symbols, so $U$ is the band-sum of the links $\{ E_k\}$ and $W$.

To summarize, after adding the additional factor $W$ to ${\gamma}'$ (and keeping the same notation, ${\gamma}'$, for the resulting curve), we have achieved $[{\gamma}]\, =\, [{\gamma}']$ in the free group $F_4$. Therefore the curves ${\gamma}, {\gamma}'$  are homotopic in the complement to the $4$-component unlink $U$. Lemma \ref{link homotopy lemma} shows that this homotopy may be assumed to be of a special kind, useful in the conclusion of the proof of the theorem further below.

\begin{lemma} \label{link homotopy lemma} {\sl
Suppose two $n$-component links $L_0, L_1$ are link-homotopic. Then there is a level-preserving link-homotopy joining them, $H\co \coprod^n S^1 \times [0,1]\longrightarrow S^3\times [0,1]$,  such that: 

1) $H$ is an isotopy  of $L_0$ for $t\in[0,\frac{1}{2})$, and \\
2) $H$ is an isotopy of $L_1$ for $t\in(\frac{1}{2},1]$ which is supported in an arbitrarily small neighborhood of a $1$-complex.

In particular, all singularities (self-intersections of link components) take place at time $t=1/2$.}
\end{lemma}

\begin{figure}[ht]
%\centering
%\vspace{.5cm} 
\includegraphics[width=14cm]{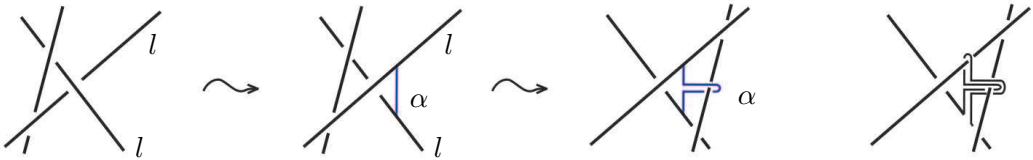}
\small
%\put(-414,2){$L_0$}
%\put(-299,2){$L_1$}
\put(-346,5){$l$}
\put(-233,5){$l$}
\put(-341,44){$l$}
\put(-229,44){$l$}
\put(-242,24){${\alpha}$}
\put(-118,25){${\alpha}$}
\caption{The two arrows represent a link homotopy followed by an isotopy. The curve ${\alpha}$ records the link homotopy data. The right-most figure in the panel shows the isotopy of $L_0$ used in Lemma \ref{link homotopy lemma}.}
\label{fig:Link homotopy}
\end{figure}

{\em Proof of lemma \ref{link homotopy lemma}.} Start with a level-preserving link homotopy $H'$ between $L_0$ and $L_1$. Each double point is formed by a transverse self-intersection of a link component. Given a double point of a link component $l$ at time $t_i$, consider an arc ${\alpha}$ in $S^3\times \{t_0+{\epsilon} \}$ with both endpoints on $l$ (and with interior disjoint from the link), such that contracting ${\alpha}$ gives the singular link at time $t_i$, see figure \ref{fig:Link homotopy}. 

At times $t>t_i+{\epsilon}$ let $\alpha$ evolve by an isotopy in the complement of the link and in the complement of other arcs. (The singular set of a link homotopy, consisting of self-intersections of components, is zero-dimensional, and ${\alpha}$ is disjoint from it by general position.) The end result at time $1$ is the link $L_1$ with a collection of arcs $\{ {\alpha} \}$, one for each double point of the link homotopy $H'$. Now a link isotopic to $L_0$ is obtained by taking a component $l$ near one of the endpoints of each curve ${\alpha}$ and performing a finger move along ${\alpha}$, as shown on the right in figure \ref{fig:Link homotopy}. This describes an isotopy of $L_0$ for $t\in[0,\frac{1}{2})$. The self-intersections and the isotopy of $L_1$ for $t\in(\frac{1}{2},1]$ consist of the reverse of the above finger moves; they are supported in a neighborhood of the arcs ${\alpha}$.
\qed

\begin{figure}[ht]
%\centering
%\vspace{.5cm} 
\includegraphics[height=9cm]{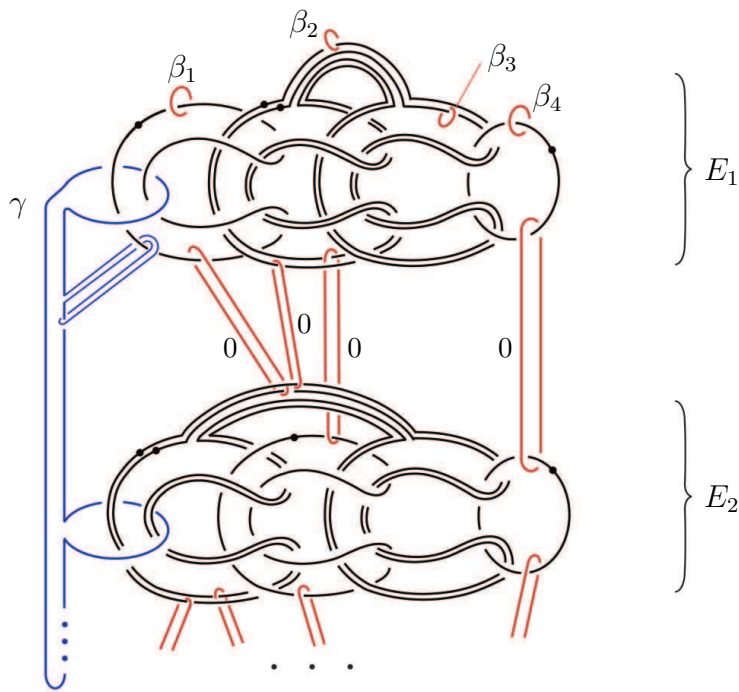}
\put(3,197){$E_1$}
\put(3,73){$E_{2}$}
\put(-260,185){${\gamma}$}
\put(-200,236){${\beta}_{1}$}
\put(-154,253){${\beta}_{2}$}
\put(-79,240){${\beta}_{3}$}
\put(-62,225){${\beta}_{4}$}
\small{
\put(-179,130){$0$}
\put(-151,139){$0$}
\put(-132,130){$0$}
\put(-75,130){$0$}
}
\caption{The band sums in figure \ref{fig:EngelProduct1} are replaced with canceling $1$-, $2$-handle pairs.}
\label{fig:EngelProduct2}
\end{figure}

Returning to the proof of the theorem, apply lemma \ref{link homotopy lemma} to the link-homotopy from ${\gamma}\cup U$ to ${\gamma}'\cup U$. Here the components of $U$ move by an isotopy, and only the curve $\gamma$ undergoes self-intersections. Consider the curve produced in the proof of lemma \ref{link homotopy lemma} at time $\frac{1}{2}-{\epsilon}$: it is isotopic to $\gamma$, and it differs from ${\gamma}'$ by finger moves along arcs $\{ \alpha\}$. Applying an isotopy to them if necessary, these arcs (and the resulting finger moves) may be assumed to be disjoint from the bands used in band-sums of the links in figure \ref{fig:EngelProduct1}. The upshot is that ${\gamma}$ is isotopic to a curve respecting the band sum decomposition, like ${\gamma}'$ in figure \ref{fig:EngelProduct1}, except for the finger moves schematically shown in figure \ref{fig:EngelProduct2}.

\begin{figure}[ht]
%\centering
%\vspace{.5cm} 
\includegraphics[height=9cm]{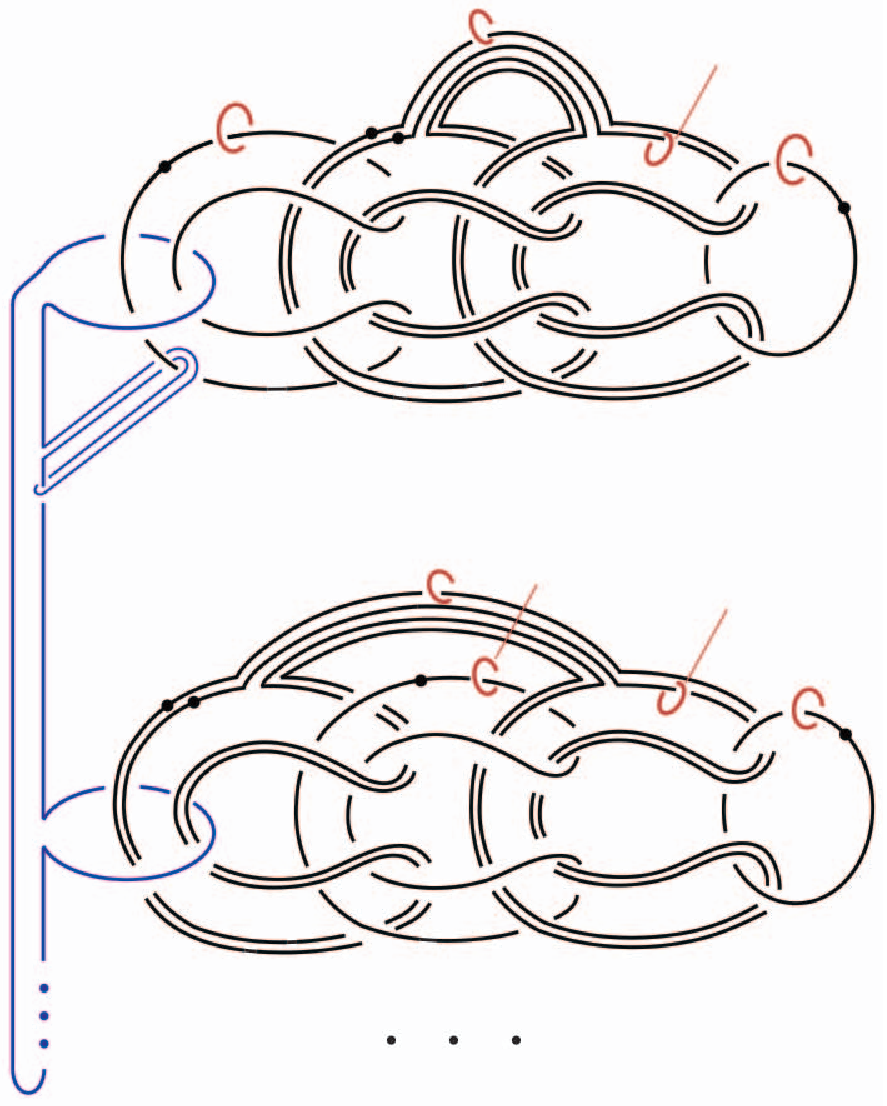}
\put(-214,185){${\gamma}$}
{
\put(-165,236){${\beta}^1_{1}$}
\put(-112,253){${\beta}^1_{2}$}
\put(-38,240){${\beta}^1_{3}$}
\put(-16,225){${\beta}^1_{4}$}
\put(-118,129){${\beta}^{2}_{1}$}
\put(-34,116){${\beta}^{2}_{3}$}
\put(-78,123){${\beta}^{2}_{2}$}
\put(-15,100){${\beta}^{2}_{4}$}
}
\caption{The zero-framed $2$-handles in figure \ref{fig:EngelProduct2}, resulting from band sums, are traded for linking circles ${\beta}^i_j$ to the components of $E_i$, bounding $2$-stage capped gropes.}
\label{fig:EngelProduct3}
\end{figure}

Now change the handle structure (but keep the $4$-manifold intact)  in figure \ref{fig:EngelProduct1} by replacing the band sums with canceling $1$-, $2$-handle pairs, as shown in figure \ref{fig:EngelProduct2}.

For the next step in the construction recall that the components of the top elementary Engel link $E_1$ have linking circles, 
the attaching curves ${\beta}_i$  of the higher stages of the capped grope B. Slide the zero-framed $2$-handles, connecting the components of $E_1$ and $E_2$, off of $E_1$, using the curves ${\beta}_i$. Now the component of $E_2$ have linking circles, denoted ${\beta}^2_i$, which bound  parallel copies of the capped gropes attached to ${\beta}_i$.  Continue this construction by sliding the zero-framed $2$-handles connecting the links $E_2$ and $E_3$ using ${\beta}^2_i$, etc. The result is that the components of each link $E_k$ have dual circles ${\beta}^k_i$, bounding a $2$-stage capped grope. (Note that the bodies of all these capped gropes are disjoint, but there are cap-cap intersections.) Here the curves ${\beta}^1_i$ are defined to be the original ${\beta}_i$. This is depicted in figure \ref{fig:EngelProduct3}. The last step of this construction yields dual circles ${\beta}^W_i$, bounding $2$-stage capped gropes, for the components of the bottom link $W$ corresponding to the factor $W$ in equation (\ref{gamma}).

%These dual capped gropes can also be arranged to be disjoint, using contraction/push-off \cite[Section 2.3]{FQ}.

%\begin{figure}[ht]
%\centering
%\vspace{.5cm} 
%\includegraphics[width=10cm]{representative.eps}
%\caption{}
%\label{fig:representative}
%\end{figure}

Now replace all dotted components  in figure \ref{fig:EngelProduct3} with zero-framed $2$-handles. The boundary $3$-manifold stays the same, but this creates new hyperbolic pairs. To solve the original surgery problem at the start of the proof of the theorem, it would suffice to solve the new stabilized surgery problem - thus finding a $4$-manifold which is a homotopy $1$-complex, and with the required boundary $3$-manifold. 

\begin{figure}[ht]
%\centering
%\vspace{.5cm} 
\includegraphics[width=7.5cm]{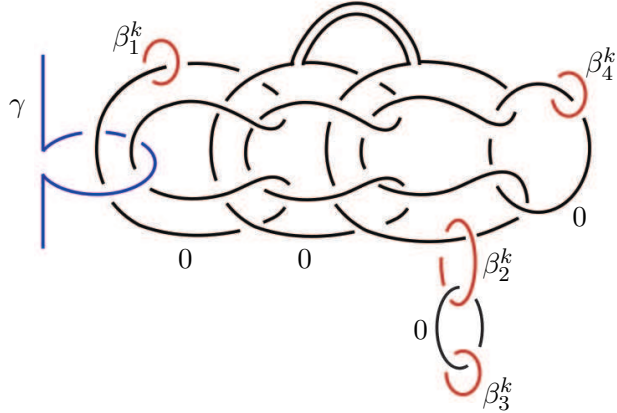}
\small
\put(-224,110){${\gamma}$}
\put(-185,135){${\beta}_1^k$}
\put(-45,49){${\beta}_2^k$}
\put(-45,0){${\beta}_3^k$}
\put(-6,125){${\beta}_4^k$}
\put(-160,51){$0$}
\put(-115,51){$0$}
\put(-71,24){$0$}
\put(-11,67){$0$}
\caption{The result of a handle slide.}
\label{fig:slide}
\end{figure}

A key feature of the elementary Engel links is that two of the components are parallel, and removing one of them yields an h-trivial link. For each link $E_k$ perform a handle slide of one of these components along its parallel. As shown in figure \ref{fig:slide}, the resulting zero-framed components which participated in the handle slide share a dual capped grope (attached to the curve labeled ${\beta}_2^k$ in the figure). In addition, the slid component has its own geometric dual, labeled ${\beta}_3^k$. The ${\delta}_{ij}$-linking  pairing with the dual capped gropes is restored by performing the handle slide of ${\beta}_2^k$ over ${\beta}_3^k$.
The result is a split link: an h-trivial link $\overline E_k$, obtained from an elementary Engel link $E_k$ by omitting a component, and a trivial unlinked component.

The handle-slide step is performed on each link $E_j$, but not on the link $W$.
To summarize,  the entire zero-framed link now consists of 
\begin{equation} \label{L} 
L:={\gamma}\, \cup\, \left( \cup_k \overline E_k\right) \, \cup({\rm unlinked \; \,  trivial \; \, components}) \, \cup\,  W,\end{equation} 
and each component of $L$ has a dual circle which is the attaching region of a $2$-stage capped grope. (For ${\gamma}$ this capped grope is given by the side $A$ of the original surgery kernel, for all other components they are parallel copies of the top two stages of $B$.) 

The main point of the band-sum and handle-slide construction above is that the link $L$ is h-trivial. Each link $\overline E_k$ is h-trivial by construction, and $W$ is h-trivial by definition, since it corresponds to a product of Milnor's relations.

Consider the disjoint null-homotopies $\Delta$ for $L$ in $D^4$. These null-homotopies, capped off with the cores of the zero-framed $2$-handles attached to $D^4$ along $L$, form the ${\pi}_1$-null spheres in the $1/2$-${\pi}_1$-null models (figure \ref{fig:NewUniversal}). Each of these spheres has a dual capped grope which constitutes the other side of the surgery models. 

The double point loops of $\Delta$ are null-homotopic in the $4$-ball; these null-homotopies intersect ${\Delta}$ and each other. All these intersections are eliminated using the dual capped gropes; now the double point loops of the spheres bound capped gropes as required in the definition of $1/2$-${\pi}_1$-null models. 	All of the gropes can be made disjoint by removing the cap-cap intersections using the usual operations of grope height raising and contraction/push-off \cite[Sections 2.3, 2.7]{FQ}.
\qed

\section{The ${\pi}_1$-Null Disk Lemma} \label{NDL section}

In this section we show that the methods developed in the proof of theorem \ref{Main theorem} also give a new formulation of the  ${\pi}_1$-null disk lemma. This lemma underlies the proofs of the surgery and s-cobordism conjecture for {\em good} groups;  indeed, the class of good groups is defined precisely as the groups for which this lemma holds (see below). Recall its statement \cite{FT}:

\smallskip

{\bf $\mathbf{{\pi}_1}$-null disk lemma (NDL)}. {\sl Let $G^c$ be a Capped Grope of height $2$ with the attaching curve $\gamma$, and ${\phi}\co {\pi}_1 G^c\longrightarrow {\pi}$ a group homomorphism. Then $\gamma$ bounds a disk ${\Delta}\co D^2\longrightarrow G^c$ which is ${\pi}_1$-null under $\phi$.}

\smallskip

A group ${\pi}$ is called {\em good} if it satisfies NDL for all choices of capped gropes $G^c$ and homomorphisms ${\phi}$.
In this statement, following the notation of \cite{FT}, capital letters in the term Capped Grope indicate the standard $4$-dimensional thickening of the underlying $2$-complex. 
Before stating the new version of NDL, consider the following rephrasing of its statement. Let $H$ denote the body of a capped grope $G^c$ of height $2$. Then $H$ is a handlebody $\natural_{i=1}^k (S^1\times D^3)$, and there is a collection of framed curves in $\partial H$ which serve as attaching regions of kinky handles (thickenings of caps of $G^c$). Connected to a basepoint in $H$, these curves represent generators of the free group $F_k\cong {\pi}_1 (H)$. Then a basic calculation in Milnor's theory of link homotopy shows that the attaching curve $\gamma$ of $G^c$ is non-trivial in the free Milnor group $M{\pi}_1(H)$. In other words, representing $H$ by  the Kirby diagram $U$ consisting of $k$ dotted components corresponding to the chosen generators of ${\pi}_1(H)$,   ${\gamma}\cup U$ is an h-essential link.  The example with all surface stages of $G^c$ of genus $1$ is shown in figure \ref{fig:wedge1}, where in the present context the curves ${\beta}_i$ serve as the attaching curves of caps of $G^c$.

Our new formulation shows that it suffices to consider curves  in the boundary  of the handlebody $H$ which are {\em trivial} in $M{\pi}_1(H)$. Consider the following statement:

\smallskip

{\bf Weak $\mathbf{{\pi}_1}$-Null disk Lemma (WNDL)}. {\sl Let $M=H \cup \, ({\rm plumbed}\; \, 2{\rm-handles})$, where $H$ is a handlebody $\natural_{i=1}^k (S^1\times D^3)$, and the $2$-handles are attached along a collection of standard framed curves in $\partial H$ which represent a set of generators of the free group ${\pi}_1(H)$. Let $\gamma$ be a curve in $\partial H$, trivial in $M{\pi}_1(H)$, where the Milnor group is defined with respect to the given generators. Let ${\phi}\co {\pi}_1 M \longrightarrow {\pi}$ be a group homomorphism. Then $\gamma$ bounds a disk ${\Delta}\co D^2\longrightarrow M$ which is ${\pi}_1$-null under $\phi$.}

\smallskip

\begin{figure}[ht]
%\centering
%\vspace{.5cm} 
\includegraphics[width=7.5cm]{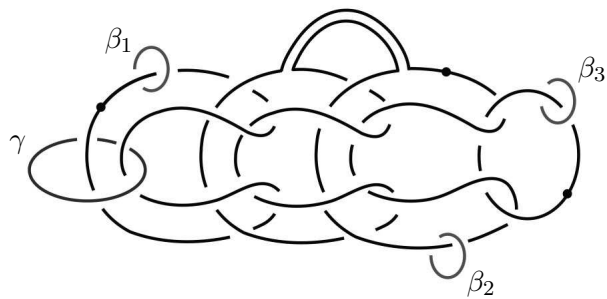}
\small
\put(-219,52){${\gamma}$}
\put(-183,90){${\beta}_1$}
\put(-46,-2){${\beta}_2$}
\put(-4,79){${\beta}_3$}
\caption{An example of ${\gamma}\subset \partial H$ in the statement of WNDL.}
\label{fig:slide1}
\end{figure}

In analogy with the ${\pi}_1$-null disk lemma, we say that a group ${\pi}$ satisfies WNDL if the statement holds for all $M, {\gamma}$ and $\phi$.
Arbitrary plumbings and self-plumbings are allowed among the $2$-handles attached to $H$ in the statement of WNDL, in analogy with caps of a capped grope.
A representative example of ${\gamma}\subset \partial (\natural_{i=1}^3 (S^1\times D^3))$  is shown in figure \ref{fig:slide1}. This example is obtained by deleting a component of an elementary Engel link, figure \ref{Elementary links figure}.

If ${\gamma}$ were null-homotopic in the complement of the dotted unlink $U$, of course it would bound a ${\pi}_1$-null disk in $M$ since $H$ is ${\pi}_1$-null in $M$. However in non-trivial cases, such as the one in figure \ref{fig:slide1}, self-intersections of the  components of $U$ are required to make ${\gamma}$ null-homotopic. Such self-intersections do not make sense since $U$ represents standard slices removed from $D^4$.

\begin{thm} \label{pi1null thm} \sl 
If a group ${\pi}$ satisfies WNDL, then ${\pi}$ is good.
\end{thm}

{\em Proof.} Consider $(G^c, {\gamma})$ and ${\phi}\co {\pi}_1 G^c\longrightarrow {\pi}$ in the statement of NDL.
The proof of theorem \ref{Main theorem} through the point (shown in figure \ref{fig:EngelProduct2}) just before stabilization applies to the curve ${\gamma}$ within the handlebody $H$. Here the circles ${\beta}_i$ linking the components, shown in figure \ref{fig:wedge1} and used throughout the proof of theorem \ref{Main theorem}, are understood to be the attaching curves of kinky handles. Instead of stabilizing (replacing dotted components with zero-framed $2$-handles), perform the handle slides on $1$-handles rather than on $2$-handles. One gets links of the form  shown in figure \ref{fig:slide}, but the components are dotted rather than $0$-framed. As in the proof of theorem \ref{Main theorem}, the link $L$ in equation (\ref{L}) is h-trivial. Since ${\pi}$ is assumed to satisfy WNDL, ${\gamma}$ bounds a ${\pi}_1$-null disk under $\phi$.
\qed

The main open case is 
${\pi}=F_k$, and ${\phi}=$ the identity homomorphism. In other words, the question is whether ${\gamma}$ as in the statement of WNDL bounds a ${\pi}_1$-null disk in $M$. For example, one may consider $k=3$ and the curve $\gamma$ in figure \ref{fig:slide}, for various choices of plumbed and self-plumbed $2$-handles attached to the curves ${\beta}_i$.


\begin{thebibliography}{10}

%[Bu02]
\bibitem{Burnside} W. Burnside, {\em On groups in which every two conjugate operations are permutable},
Proc. Lond. Math. Soc. 35 (1902), 28-37.


%[CF84]
\bibitem{CF} A. Casson and M.  Freedman, {\em Atomic surgery problems}, Four-manifold theory (Durham, N.H., 1982), 181-199, Contemp. Math., 35, Amer. Math. Soc., Providence, RI, 1984.


\bibitem{Cochran} T.D. Cochran, Derivatives of links: Milnor's concordance invariants and Massey's products. Mem. Amer. Math. Soc. 84 (1990), no. 427.

%[Fr82a]
\bibitem{F0} M. Freedman, {\em The topology of four-dimensional
manifolds}, J. Differential Geom. 17(1982), 357-453.


%[Fr83]
\bibitem{F1} M. Freedman, {\em The disk theorem for
four-dimensional manifolds}, Proc. ICM Warsaw (1983), 647-663.

\bibitem{FK} M. Freedman and V.  Krushkal, 
{\em Engel relations in 4-manifold topology},
Forum Math. Sigma 4 (2016), e22, 57 pp.

%[FQ90]
\bibitem{FQ} M. Freedman and F. Quinn, {\em The topology of
4-manifolds}, Princeton Math. Series 39, Princeton, NJ, 1990.

%[FT95a]
\bibitem{FT} M. Freedman and P. Teichner,
{\em$4$-Manifold Topology I: Subexponential groups}, Invent. Math. 122 (1995), 509-529.

\bibitem{FT1} M. Freedman and P. Teichner,
{\em$4$-Manifold Topology II:}, Invent. Math. 122 (1995)


%[Hei61]
\bibitem{Hei} H. Heineken, 
{\em Engelsche Elemente der L\'{a}nge drei}, 
Illinois J. Math. 5 (1961), 681-707. 

%[Ho29]
\bibitem{Hopkins} C. Hopkins, {\em Finite groups in which conjuate operations are commutative}, Am. J.
Math. 51, (1929), 35-41.


%[KQ00]
\bibitem{KQ} V. Krushkal and F. Quinn, {\em Subexponential groups in
$4$-manifold topology}, Geom. Topol. 4 (2000), 407-430.

%[Lev42]
\bibitem{Levi} F.W. Levi,
{\em Groups in which the commutator operation satisfies certain algebraic conditions}, 
J. Indian Math. Soc. 6 (1942), 87-97.


%[Mil54]
\bibitem{M} J. Milnor, {\em Link Groups}, Ann. Math 59 (1954), 177-195.


\bibitem{Wall} C.T.C. Wall, Surgery on compact manifolds. London Mathematical Society Monographs, No. 1. Academic Press, London-New York, 1970.




\end{thebibliography}
\end{document}